\documentclass{amsart}
\usepackage[latin1]{inputenc}
\usepackage[T1]{fontenc}
\usepackage[english]{babel}

\usepackage{amsmath,amssymb,amsthm,amsfonts}

\theoremstyle{plain}
\newtheorem{theorem}{Theorem}[section]
\newtheorem{corollary}[theorem]{Corollary}
\newtheorem{prop}[theorem]{Proposition}
\newtheorem{lemma}[theorem]{Lemma}
\theoremstyle{definition}
\newtheorem{remark}[theorem]{Remark}
\newtheorem{remarks}[theorem]{Remarks}
\newtheorem{example}[theorem]{Example}

\newtheorem{definitions}[theorem]{Definitions}

\newcommand{\D}{\mathbb{D}}
\newcommand{\R}{\mathbb{R}}

\newcommand{\N}{\mathbb{N}}
\newcommand{\T}{\mathbb{T}}

\newcommand{\eps}{\varepsilon}
\newcommand{\ext}[1][X^*]{\ensuremath{\mathrm{ext}(B_{#1})}}

\newcommand{\e}{\mathrm{e}}
\newcommand{\Id}{\mathrm{Id}}
\newcommand{\dist}{\mathrm{dist}}

\newcommand{\conv}{\mathrm{co}}
\newcommand{\aconv}{\mathrm{aco}}
\newcommand{\Lin}{\mathrm{lin}}

 \DeclareMathOperator{\re}{Re}

\renewcommand{\leq}{\leqslant}
\renewcommand{\geq}{\geqslant}

\renewcommand{\le}{\leqslant}

\newcounter{equi1}
\newenvironment{equi}
{\begin{list}
        {(\roman{equi1})}
        {\setlength{\itemsep}{0ex plus 0.2ex minus 0ex}
         \setlength{\topsep}{0ex}
         \setlength{\parsep}{0ex}
         \setlength{\labelwidth}{7ex}
         \usecounter{equi1}}}
{\end{list}}

\begin{document}
\begin{center}\small
[Illinois J.\ Math.\ (to appear)]
\end{center}

\title[Convexity and smoothness of spaces with numerical index~$1$]
{Convexity and smoothness of Banach spaces with numerical index
one}

 \date{July 8th, 2008. Revised form: November 3rd, 2008.}

 \thanks{The first author was supported by Junta de
Andaluc\'{\i}a grant P06-FQM-01438. The second and the fourth
authors were partially supported by the Spanish MEC project
no.\ MTM2006-04837 and Junta de Andaluc\'{\i}a grants FQM-185 and
P06-FQM-01438. The third author was partially supported by a
Juan de la Cierva grant, by the Spanish MEC project no.\
MTM2006-04837 and Junta de Andaluc\'{\i}a grants FQM-185 and
P06-FQM-01438.}
\author[Kadets]{Vladimir Kadets}
\author[Mart\'{\i}n]{Miguel Mart\'{\i}n}
\author[Mer\'{\i}]{Javier Mer\'{\i}}
\author[Pay\'{a}]{Rafael Pay\'{a}}

\address[Kadets]{Department of Mechanics and Mathematics,
Kharkov National University, pl.~Svobody~4,  61077~Kharkov,
Ukraine} \email{\texttt{vova1kadets@yahoo.com}}

\address[Mart\'{\i}n, Mer\'{\i}, Pay\'{a}]{Departamento de An\'{a}lisis Matem\'{a}tico \\ Facultad de
 Ciencias \\ Universidad de Granada \\ 18071 Granada, Spain}
\email{\texttt{mmartins@ugr.es}, \texttt{jmeri@ugr.es},
\texttt{rpaya@ugr.es} }

\maketitle

\begin{abstract}
We show that a Banach space with numerical index one cannot
enjoy good convexity or smoothness properties unless it is
one-dimensional. For instance, it has no WLUR points in its
unit ball, its norm is not Fr\'{e}chet smooth and its dual norm is
neither smooth nor strictly convex. Actually, these results
also hold if the space has the (strictly weaker) alternative
Daugavet property. We construct a (non-complete) strictly
convex predual of an infinite-dimensional $L_1$ space (which
satisfies a property called lushness which implies numerical
index~$1$). On the other hand, we show that a lush real Banach
space is neither strictly convex nor smooth, unless it is
one-dimensional. Therefore, a rich subspace of the real space
$C[0,1]$ is neither strictly convex nor smooth. In particular,
if a subspace $X$ of the real space $C[0,1]$ is smooth or
strictly convex, then $C[0,1]/X$ contains a copy of $C[0,1]$.
Finally, we prove that the dual of any lush
infinite-dimensional real space contains a copy of $\ell_1$.
\end{abstract}

\section{Introduction}
The classical formula $ \|T\| = \sup\{\left|\langle Tx, x
\rangle \right|\,:\, x \in X,\ \|x\|=1\} $ for the norm of a
self-adjoint operator $T$ on a Hilbert space $X$ can be
rewritten, thanks to the well-known representation of the dual
$X^*$ as
\begin{equation} \label{mur-1-eq1}
\|T\| = \sup\{|x^*(Tx)|\,:\, x \in X,\ x^* \in X^*,\ x^*(x)=\|x^*\|=\|x\|=1\}.
\end{equation}
For a non self-adjoint operator this formula may fail.
Nevertheless, there are some Banach spaces $X$ in which
equality \eqref{mur-1-eq1} is valid for \emph{every} bounded
linear operator $T$ on $X$. As we will explain below, such
spaces are said to have numerical index~$1$. Among these spaces
are all classical $C(K)$ and $L_1(\mu)$ spaces.

Given a real or complex Banach space $X$, we write $B_X$,
$S_X$, $X^*$ and $L(X)$, to denote, respectively, the closed
unit ball, the unit sphere, the topological dual and the Banach
algebra of bounded linear operators on $X$.

The \emph{numerical range} of an operator $T\in L(X)$ is the
subset of the base field given by
$$
V(T)=\{x^*(Tx)\ : \ x^*\in S_{X^*},\ x\in S_X,\ x^*(x)=1\},
$$
and the \emph{numerical radius} of $T$ is then given by
$v(T)=\sup\{|\lambda|\ : \ \lambda\in V(T)\}$. These concepts
were independently introduced by F.~Bauer \cite{Bauer} and
G.~Lumer \cite{Lumer} in the 1960's to extend the classical
field of values of matrices (O.~Toeplitz, 1918 \cite{Toe}). We
refer the reader to the monographs by F.~Bonsall and J.~Duncan
\cite{B-D1,B-D2} for a detailed account. The \emph{numerical
index} of the space $X$ (Lumer, 1968 \cite{D-Mc-P-W}) is the
constant $n(X)$ defined by
$$
n(X):=\inf\{v(T)\ : \ T\in L(X),\ \|T\|=1\}
$$
or, equivalently, the greatest constant $k\geq 0$ such that
$k\|T\|\leq v(T)$ for every $T\in L(X)$. Observe that
$0\leqslant n(X) \leqslant 1$ for every Banach space $X$, and
$n(X)=1$ if and only if equality \eqref{mur-1-eq1} is valid for
all operators on $X$. The reader will find the state-of-the-art
on numerical indices in the recent survey paper \cite{KaMaPa}
to which we refer for background.

Let us mention here several facts concerning the numerical
index which are relevant to our discussion. Examples of Banach
spaces having numerical index $1$ are $C(K)$ spaces, $L_1(\mu)$
spaces, Lindenstrauss spaces (i.e.\ isometric preduals of
$L_1(\mu)$ spaces) \cite{D-Mc-P-W}, all function algebras
\cite{WerJFA} (for instance, the disk algebra $A(\D)$ and
$H^\infty$), and finite-codimensional subspaces of $C[0,1]$
\cite{BKMW}. Next, one has $v(T^*)=v(T)$ for every $T\in L(X)$,
where $T^*$ is the adjoint operator of $T$ (see \cite[\S
9]{B-D1}), and it clearly follows that $n(X^*)\leqslant n(X)$
for every Banach space $X$. It has recently been discovered
that this inequality can be strict. Actually, in
\cite[Example~3.1]{BKMW} an example is given of a real Banach
space $X$ such that $n(X)=1$ while $n(X^*)=0$. We refer to the
very recent paper \cite{Mart-L} for sufficient conditions to
ensure the equality in the inequality $n(X^*)\leqslant n(X)$.
Every separable Banach space containing $c_0$ can be
equivalently renormed to have numerical index~$1$
\cite[\S4]{BKMM-lush}, in particular, this happens with any
closed subspace of $c_0$. On the other hand, there is no
infinite-dimensional real reflexive space with numerical
index~$1$ \cite{LMP1999}.

Our main goal in this paper is to study which convexity or
smoothness properties are possible for the unit ball of a
Banach space with numerical index~$1$. At the end of this
introduction we give the necessary definitions of the convexity
and smoothness properties we use along the paper. A difficulty
with such a study is that the property of having numerical
index $1$ deals with all operators on the space and we do not
know of any characterization of it in terms of the space and
its successive duals. The previous solutions to this difficulty
have been to deal with either weaker or stronger geometrical
properties. Let us briefly give an account of some of them. Let
$X$ be a real or complex Banach space.
\begin{enumerate}
\item[$(a)$] $X$ is said to be a \emph{CL-space} if $B_X$
    is the absolutely convex hull of every maximal convex
    subset of $S_X$.
\item[$(b)$] We say that $X$ is an \emph{almost-CL-space}
    if $B_X$ is the closed absolutely convex hull of every
    maximal convex subset of $S_X$.
\item[$(c)$] $X$ is \emph{lush} if for every $x,y \in S_X$
    and every $\eps > 0$, there is a slice
    $$
    S=S(x^*,\eps):=\{z\in B_X\ : \ \re x^*(z)>1-\eps\}
    $$
    with $x^* \in S_{X^*}$ such that $x \in S$ and $\dist
    \left(y,\aconv (S)\right) < \eps$, where $\aconv(S)$
denotes the absolutely convex hull of the set $S$.
\item[$(d)$] $X$ has \emph{numerical index $1$} ($n(X)=1$
    in short) if $v(T)=\|T\|$ for every $T\in L(X)$.
\item[$(e)$] We say that $X$ has the \emph{alternative
    Daugavet property} provided that every rank-one
    operator $T\in L(X)$ satisfies $v(T)=\|T\|$. The same
    equality is then satisfied by all weakly compact
    operators on $X$ \cite[Theorem~2.2]{MaOi}.
\end{enumerate}
The implications $(a)\Longrightarrow (b) \Longrightarrow (c)$
and $(d)\Longrightarrow (e)$ are clear and none of them
reverses (see \cite[\S3 and \S7]{BKMW} for a detailed account).
Also, $(c)\Longrightarrow (d)$ by \cite[Proposition~2.2]{BKMW}.

Some additional comments on the above properties may be in
place. CL-spaces where introduced in 1960 by R.~Fullerton
\cite{Full} and it was later shown that a finite-dimensional
Banach space has numerical index~$1$ if and only if it is a
CL-space (\cite[Theorem~3.1]{Mc} and
\cite[Corollary~3.7]{Lima2}). Therefore, the above five
properties are equivalent in the finite-dimensional case. All
$C(K)$ spaces as well as real $L_1(\mu)$ spaces are CL-spaces,
while infinite-dimensional complex $L_1(\mu)$ spaces are only
almost-CL-spaces (see \cite{MartPaya-CL}). Almost-CL-spaces
first appeared without a name in the memoir by J.~Lindenstrauss
\cite{Lin-book} and were further discussed by \AA.~Lima
\cite{Lima,Lima2} who showed that real Lindenstrauss spaces
(i.e.\ isometric preduals of $L_1(\mu)$) are CL-spaces
\cite[\S3]{Lima} and complex Lindenstrauss spaces are
almost-CL-spaces \cite[\S3]{Lima2}. The disk algebra is another
classical example of an almost-CL-space
\cite[Theorem~32.9]{B-D2}. More information can be found in
\cite{Cheng-Li,MarRNP,MartPaya-CL,Reis}.

Lush spaces were introduced recently \cite{BKMW} and they were
the key to provide an example of a Banach space $X$ such that
$n(X^*)<n(X)$ and to estimate the polynomial numerical index of
some spaces \cite{ChoGarMaeMar,KimMartMeri}. We refer to
\cite{BKMM-lush} for characterizations and examples of lush
spaces. Among the advantages of the concept of lushness are
that this property is separably determined
\cite[Theorem~4.2]{BKMM-lush} and that it gives many new
examples of Banach spaces with numerical index~$1$. Namely,
C-rich subspaces of $C(K)$ are lush and so they have numerical
index~$1$ \cite[Theorem~2.4]{BKMW}. A closed subspace $X$ of a
$C(K)$ space is said to be \emph{C-rich} if for every nonempty
open subset $U$ of $K$ and every $\eps > 0$, there is a
positive function $h\in C(K)$ of norm $1$ with support inside
$U$ such that the distance from $h$ to $X$ is less than $\eps$.
This definition covers finite-codimensional subspaces of
$C[0,1]$ \cite[Proposition~2.5]{BKMW}, so they are lush. Also,
all function algebras are lush (see
\cite[Example~2.4]{BKMM-lush} and \cite[\S3]{WerJFA}).

The alternative Daugavet property was introduced and
characterized in \cite{MaOi} but, in an equivalent way, the
property defining it had appeared in some papers of the 1990's.
The name comes from the fact that an operator $T$ on a Banach
space $X$ satisfies $v(T)=\|T\|$ if and only if $\|\Id +
\omega\,T\|=1 + \|T\|$ for some $\omega\in \T$ ($\T$ being the
set of modulus one scalars) \cite{D-Mc-P-W}, that is, the
operator $S=\omega\,T$ satisfies the so-called Daugavet
equation $\|\Id + S\|=1 + \|S\|$. Therefore, $X$ has the
alternative Daugavet property if and only if every rank-one
operator (equivalently, every weakly compact operator)
satisfies the Daugavet equation up to rotation. We refer to the
already cited paper \cite{MaOi} and to \cite{Mar-ADP} for more
information and background. Let us comment that Banach spaces
with the Radon-Nikod\'{y}m property and the alternative Daugavet
property are actually almost-CL-spaces \cite{MarRNP}. Asplund
spaces with the alternative Daugavet property are lush, but
they need not be almost-CL-spaces \cite[Example~2.4]{BKMW}.

The main question in this paper, not yet solved, is whether a
Banach space $X$ with $n(X)=1$ can be smooth or strictly
convex. Two remarks are pertinent. First, even though the exact
value of $n(\ell_p^2)$ is not known for $p\neq 1,2,\infty$ (see
\cite{MarMer-Lp}), it is known that the set
$$
\bigl\{n\bigl(\ell_p^{2}\bigr)\ : \ 1<p<\infty\bigr\}
$$
contains all possible values of the numerical index except $1$
\cite{D-Mc-P-W}. Thus, the question above only makes sense for
the value $1$. Second, it is clear that an almost-CL-space
cannot be strictly convex (almost-CL-spaces are somehow the
extremely opposite property to strict convexity), and it has
recently been shown that a real almost-CL-space cannot be
smooth \cite[Theorem~3.1]{Cheng-Li}.

Let us summarize the main results in this paper.
Section~\ref{sec:occlusive-ADP} is devoted to show that a
Banach space with the alternative Daugavet property and
dimension greater than one has no WLUR points in its unit ball,
its norm is not Fr\'{e}chet smooth and its dual norm is neither
strictly convex nor smooth. Next, in \S\ref{sec:non-complete}
we construct a non-complete predual of an $L_1(\mu)$ space
which is strictly convex. This space is lush (extending this
definition to general normed spaces literally), while its
completion is an almost-CL-space. The aim of
section~\ref{sec:lush=>stronglylush} is to show that separable
lush spaces actually satisfy a stronger property: there is a
norming subset $\tilde K$ of $S_{X^*}$ such that for every
$x^*\in \tilde{K}$ and every $\eps>0$, one has
$$
B_X=\overline{\aconv\bigl(S(x^*,\eps)\bigr)}.
$$
In the real case, it is actually true that $B_X$ is the closed
absolutely convex hull of the (non-empty) face generated by
$x^*$. This implies that a real lush Banach space is neither
strictly convex nor smooth, unless it is one-dimensional.
Therefore, a C-rich subspace of the real space $C[0,1]$ is
neither strictly convex nor smooth, and this answers a question
of M.~Popov from 1996. In particular, if a subspace $X$ of the
real space $C[0,1]$ is smooth or strictly convex, then
$C[0,1]/X$ contains a copy of $C[0,1]$. We devote
\S\ref{sec:localizations} to some localizations of convexity
and smoothness properties. Namely, it was asked in
\cite[Problem~13]{KaMaPa} whether a Banach space $X$ with
$n(X)=1$ satisfies that $|x^*(x)|=1$ for every $x\in \ext[X]$
(the set of extreme points in $B_X$) and every $x^*\in
\ext[X^*]$, as it happens in the finite-dimensional case
\cite{Mc} (a positive answer would lead to the impossibility of
having a strictly convex space with numerical index~$1$ other
than the one-dimensional one). But actually, we construct
examples of separable lush spaces where this does not happen,
giving a negative answer to the cited problem. On the other
hand, we show that for lush spaces, $|x^*(x)|=1$ for every
$x^*\in \ext[X^*]$ and every $w^*$-extreme point, which gives
us that a lush space which is WMLUR has to be one-dimensional.

We finish the introduction with the definitions and notations
of the convexity and smoothness properties that we need
throughout the paper. We refer the reader to the books
\cite{DGZ,DiesGeom} and the papers
\cite{BandHuandLinTroyans,Kunen-Rosenthal} for more information
and background.

The norm of a real or complex Banach space $X$ (or $X$ itself)
is said to be \emph{smooth} if for every $x\in S_X$, there is a
unique norm-one functional $x^*$ such that $x^*(x)=1$. The
space $X$ is said to be \emph{strictly convex} when
$\ext[X]=S_X$. It is well known that $X$ is smooth (resp.\
strictly convex) if $X^*$ is strictly convex (resp.\ smooth),
but the converse is not true. We say that the norm of $X$ is
\emph{Fr\'{e}chet smooth} when the norm of $X$ is Fr\'{e}chet
differentiable at any point of $S_X$. By the Smulyan test, the
norm of $X$ is Fr\'{e}chet smooth if and only if every functional
$x^*\in S_{X^*}$ which attains its norm is $w^*$-strongly
exposed (i.e.\ there is $x\in S_X$ such that for every sequence
$(x^*_n)$ in $B_{X^*}$ such that $x^*_n(x)\longrightarrow
1=x^*(x)$ one has $x^*_n\longrightarrow x^*$ in norm).

An $x\in S_X$ is said to be a point of \emph{local uniform
rotundity (LUR point)} if $\|x_n -x\|\longrightarrow 0$ for
every sequence $(x_n)$ in $S_X$ such that
$\|x_n+x\|\longrightarrow 2$. If for every sequence $(x_n)$ of
$S_X$ with $\|x_n+x\|\longrightarrow 2$ one only has that
$x_n\longrightarrow x$ in the weak topology, we say that $x$ is
a point of \emph{weakly local uniform rotundity (WLUR point)}.

A point $x$ in $S_X$ is said to be (\emph{weakly})
\emph{midpoint locally uniformly rotund} or \emph{MLUR} (resp.\
\emph{WMLUR}) if for any sequence $(y_n)$ in $B_X$, $\lim_n \|x
\pm y_n\| \leq 1$ implies $\lim_n \|y_n\| = 0$ ($\lim_n y_n =
0$ in the weak topology). A point $x$ of $B_X$ is called
\emph{weak$^*$-extreme} if it is an extreme point of
$B_{X^{**}}$. Every WMLUR point of $B_X$ is a weak$^*$-extreme
point of $B_X$ (see \cite[p.~674]{GodLinTroy} and
\cite[p.~173]{Kunen-Rosenthal}). We say that the norm of $X$ is
MLUR (WMLUR) if every point in $S_X$ is MLUR (WMLUR).

\section{Prohibitive results for the alternative Daugavet
property}\label{sec:occlusive-ADP} The aim in this section is
to show that there are some convexity and smoothness properties
which are incompatible with the alternative Daugavet property
and so, they are incompatible with the numerical index~$1$. We
start with smoothness and strict convexity of the dual norm.

\begin{theorem}\label{thm-dual-smooth-strictly-convex}
Let $X$ be a Banach space with the alternative Daugavet
property and dimension greater than one. Then, $X^*$ is neither
smooth nor strictly convex.
\end{theorem}

\begin{proof}
Since the dimension of $X$ is greater than $1$, we may find
$x_0\in S_X$ and $x_0^*\in S_{X^*}$ such that $x^*_0(x_0)=0$.
Then, we consider the norm-one operator $T=x_0^*\otimes x_0$,
which satisfies $T^2=0$. On the other hand, thanks to
\cite[Theorem~1.2]{AcoPay93}, there is a sequence of norm-one
operators $(T_n)$ converging in norm to $T$ and such that the
adjoint of each of them attains its numerical radius. Moreover,
we may suppose that all the $T_n$'s are compact by
\cite[Remark~1.3]{AcoPay93}. Since $X$ has the alternative
Daugavet property, we get
$$
v(T_n^*)=v(T_n)=\|T_n\|=1.
$$
As the operators $T_n^*$ attain their numerical radius, for
every positive integer $n$, we may find $\lambda_n\in \T$ and
$(x_n^*,x_n^{**})\in S_{X^*}\times S_{X^{**}}$ such that
\begin{equation}\label{eq:proof-dual-smooth}
\lambda_n\,x_n^{**}(x_n^*)=1 \qquad \text{ and } \qquad
\bigl[T_n^{**}(x_n^{**})\bigr](x_n^*)=x_n^{**}(T_n^*(x_n^*))=1.
\end{equation}

If $X^*$ is smooth, we deduce that
$$
T_n^{**}(x_n^{**})=\lambda_n\,x_n^{**} \qquad \bigl(n\in \N\bigr).
$$
Thus,
$$
\left\|\left[T_n^{**}\right]^2(x_n^{**}) \right\| =\|\lambda_n^2\,x_n^{**}\|=1 \qquad (n\in \N).
$$
But, since $T_n\longrightarrow T$ and $T^2=0$, we have that
$\left[T_n^{**}\right]^2 \longrightarrow 0$, a contradiction.

If $X^*$ is strictly convex, we deduce from
\eqref{eq:proof-dual-smooth} that
$$
T_n^{*}(x_n^{*})=\lambda_n\,x_n^{*} \qquad \bigl(n\in \N\bigr),
$$
which leads to a contradiction the same way as before.
\end{proof}

As a consequence of the above result, we get that $n(H^1)<1$,
where $H^1$ represents the Hardy space. Actually, we have more.

\begin{example}
{\slshape Let $X$ be $C(\T)/A(\D)$. Then, its dual $X^*=H^1$ is
smooth (see \cite[Remark IV.1.17]{HWW}, for instance), so $X$
does not have the alternative Daugavet property by
Theorem~\ref{thm-dual-smooth-strictly-convex} and neither does
$X^*=H^1$. In particular, $n(X)<1$ and $n(X^*)<1$.}
\end{example}

\begin{remarks}$ $
\begin{itemize}
\item[(a)] The proof of
    Theorem~\ref{thm-dual-smooth-strictly-convex} can be
    adapted to yield the following result. {\slshape Let
    $X$ be a Banach space with the alternative Daugavet
    property and such that the set of compact operators
    attaining its numerical radius is dense in the space of
    all compact operators. Then, $X$ is neither strictly
    convex nor smooth, unless it is one-dimensional.\ }
    Indeed, we may follow the proof of
    Theorem~\ref{thm-dual-smooth-strictly-convex} (without
    considering adjoint operators) to get the result.
\item[(b)] It is known that for Banach spaces with the
    Radon-Nikod\'{y}m property, the set of compact operators
    attaining their numerical radius is dense in the space
    of all compact operators \cite[Theorem~2.4]{AcoPay93}.
    Therefore, we get that {\slshape a Banach space having
    the Radon-Nikod\'{y}m property and the alternative Daugavet
    property is neither smooth nor strictly convex, unless
    it is one-dimensional.\ }
\item[(c)] Actually, the above result was essentially
    known. Namely, if $X$ has the alternative Daugavet
    property and the Radon-Nikod\'{y}m property, then $X$ is an
    almost-CL-space \cite[Theorem~1]{MarRNP}. It is clear
    that a (non-trivial) almost-CL-space cannot be strictly
    convex. On the other hand, the fact that a non-trivial
    real almost-CL-space cannot be smooth follows from a
    very recent result \cite[Theorem~3.1]{Cheng-Li}.
\item[(d)] The fact that there are Banach spaces in which
    the set of numerical radius attaining operators is not
    dense in the space of all operators was discovered in
    1992 \cite{Paya-counterexample}. Nevertheless,
    {\slshape we do not know of any Banach space for which
    the set of compact operators which attain their
    numerical radius is not dense in the space of all
    compact operators.}
\item[(e)] Let us comment that it is also an open problem
    whether a Banach space with the Daugavet property can
    be smooth or strictly convex. We recall that a Banach
    space has the \emph{Daugavet property} if $\|\Id +
    T\|=1 + \|T\|$ for every rank-one operator $T\in L(X)$
    \cite{KSSW}. It is clear that the Daugavet property
    implies the alternative Daugavet property (and the
    converse result is not true). Therefore, an example of
    a smooth or strictly convex Banach space with the
    Daugavet property would give an example of a Banach
    space where the rank-one operators cannot be
    approximated by compact operators attaining the
    numerical radius.
\end{itemize}
\end{remarks}

More prohibitive results for the alternative Daugavet property
are the following.

\begin{prop}
Let $X$ be a Banach space with the alternative Daugavet
property. Then, $B_X$ fails to contain a WLUR point, unless $X$
is one-dimensional.
\end{prop}

\begin{proof}
Let $x_0\in S_X$ be a WLUR point. If the dimension of $X$ is
greater than $1$, there is $x_0^*\in S_{X^*}$ such that
$x_0^*(x_0)=0$. Then, the rank-one operator $T=x_0^*\otimes
x_0$ satisfies $\|T\|=1$ and so, $v(T)=1$. Therefore, we may
find sequences $(x_n)$ in $S_X$ and $(x_n^*)$ in $S_{X^*}$ such
that
$$
x_n^*(x_n)=1 \qquad \text{ and } \qquad |x_n^*(x_0)|\,|x_0^*(x_n)|
=|x_n^*(Tx_n)|\longrightarrow 1.
$$
Therefore, we get $|x_n^*(x_0)|\longrightarrow 1$ and
$|x_0^*(x_n)| \longrightarrow 1$. If for every $n\in \N$ we
take $\lambda_n\in \T$ such that
$x_n^*(x_0)=\lambda_n\,|x_n^*(x_0)|$, we have
\begin{equation*}
2\geq \|x_0 + \lambda_n\,x_n\|\geq |x_n^*(x_0 + \lambda_n\,x_n)|\longrightarrow 2.
\end{equation*}
So, being $x_0$ a WLUR point, we get that
$(\lambda_n\,x_n)\longrightarrow x_0$ in the weak topology,
which contradicts the fact that $|x_0^*(x_n)| \longrightarrow
1$.
\end{proof}

The above result is not true if we replace the WLUR point by a
point of Fr\'{e}chet smoothness. For instance, $n(c_0)=1$ but the
norm of $c_0$ is Fr\'{e}chet differentiable at a dense subset of
$S_{c_0}$ since $c_0$ is Asplund. But it is not difficult to
show that a Banach space with the alternative Daugavet property
cannot have a Fr\'{e}chet smooth norm, unless it is
one-dimensional.

\begin{prop}
Let $X$ be a Banach space with the alternative Daugavet
property. Then, the norm of $X$ is not Fr\'{e}chet smooth, unless
$X$ is one-dimensional.
\end{prop}

\begin{proof}
Using \cite[Lemma~1]{LMP1999}, we have that $|x^{**}(x^*)|=1$
for every $x^{**}\in \ext[X^{**}]$ and every $w^*$-strongly
exposed point $x^*$ of $B_{X^*}$. Now, if the norm of $X$ is
Fr\'{e}chet-smooth, then every functional on $S_{X^*}$ attaining
its norm is $w^*$-strongly exposed (see
\cite[Corollary~I.1.5]{DGZ} for instance). Since, by the
Bishop-Phelps Theorem, we have that the set of norm-attaining
norm-one functionals is (norm) dense on $S_{X^*}$, we get that
$$
|x^{**}(x^*)|=1
$$
for all $x^{**}\in \ext[X^{**}]$ and all $x^*\in S_{X^*}$. This
clearly leads to the fact that $X$ is one-dimensional.
\end{proof}

\section{A noncomplete strictly convex lush space}\label{sec:non-complete}
The aim of this section is to construct an example of a
non-complete infinite-dimensional strictly convex normed space
with numerical index~$1$ (actually lush). As we will see, its
completion is very far away from being strictly convex. In the
next section, we will show that actually no real lush complete
space can be strictly convex.

We need some definitions and preliminary results.

\begin{definitions}
Let $|||\cdot|||$ and $\|\cdot\|$ be two norms on a linear
space $X$ and $\eps >0$. We say that $|||\cdot|||$ is
\emph{$\eps$-equivalent} to $\|\cdot\|$ if
\begin{equation*}
\frac{1}{1+\eps} \|x\| \leq |||x||| \leq (1+\eps) \|x\| \qquad (x\in X).
\end{equation*}
A property $\mathcal P$ of normed spaces is said to be a
\emph{stable C-property}, if $C[0,1] \in \mathcal P$ and for
every Banach space $X$ the following  condition is sufficient
for $X \in \mathcal P$:\newline for every $\eps > 0$ and for
every finite subset $F \subset X$, there is a subspace $Y
\subset X$, such that $F \subset Y$ and $Y$ possesses an
$\eps$-equivalent norm $\|\cdot\|_\eps$ with $\left(Y,
\|\cdot\|_\eps \right) \in \mathcal P$.
\end{definitions}

It is immediate that lushness and the alternative Daugavet
property are stable C-properties.

We are now ready to state the main result to get the example.

\begin{theorem} \label{str-conv-thm1}
For every strictly convex separable Banach space $Y_0$, there
is a strictly convex separable normed space $X \supset Y_0$
possessing all stable C-properties.
\end{theorem}

We will need the following surely well-known lemma.

\begin{lemma} \label{str-conv-lem1}
Let $Y$ be a strictly convex closed subspace of a separable
Banach space $X$. Then for every $\eps >0$, there is an
$\eps$-equivalent strictly convex norm $|||\cdot|||$ on $X$
which coincides with the original one on $Y$.
\end{lemma}

\begin{proof}
The existence of a norm $p$ satisfying all conditions of this
statement except being $\eps$-equivalent to the original one is
well known (see for example \cite[p.~84]{DGZ} or
\cite[Theorem~1.1]{Tang}). Then, for sufficiently small $t >
0$, the norm $|||x|||:= (1-t)\|x\| + tp(x)$ will be the one
which we need.
\end{proof}

\begin{proof}[Proof of Theorem~\ref{str-conv-thm1}]
We are going to construct a sequence of separable strictly
convex Banach spaces $(X_n)$ with the following properties:
\begin{equi}
\item $X_1 = Y_0$.
\item $X_n$ is a subspace of $X_m$ for $n<m$.
\item For every $n \in \N$ there is a $\frac1n$-equivalent
    norm $\|\cdot\|_{n}$ on $X_{n}$ with
    $\left(X_{n}, \|\cdot\|_{n} \right)$ being isometric to $C[0,1]$.
 \end{equi}
Since $X_1$ is already known, the only thing we need for this
construction is to show how to get $X_{m+1}$ from $X_m$. Let us
fix $m\in \N$. Since $X_{m}$ is separable, we can (and do so)
consider $X_{m}$ as a subspace of $C[0,1]$. According to
Lemma~\ref{str-conv-lem1}, there is an $(\frac1m)$-equivalent
strictly convex norm $|||\cdot|||$ on $C[0,1]$ which coincides
with the original norm on $X_{m}$.  Put $X_{m+1} =
\bigl(C[0,1], |||\cdot|||\bigr)$, and the original norm of
$C[0,1]$ plays the role of $\|\cdot\|_{n}$ in the condition
(iii). So the construction is completed.

What remains to complete the proof itself is to put $X =
\bigcup_{m\in \N} X_m$. Then, for every $\eps> 0$ and for every
finite subset $F \subset X$, one can find $n\in \N$ such that
$\frac1n < \eps$ and $F \subset X_n$. Since $\|\cdot\|_n$ is
$\eps$-equivalent to the norm of $X_n$, we get the requirement.
\end{proof}

Since lushness is a stable C-property, we get the desired
example.

\begin{example}
{\slshape There are normed lush spaces which are strictly
convex. }
\end{example}

We are going to show that the completions of the above examples
(which are of course also lush) are not strictly convex.
Actually, they are almost-CL-spaces.

Following Bourgain's book \cite{Bourgain-book}, we say that a
Banach space $X$ is an \emph{${\mathcal L}^\infty_{1+}$-space}
if for any finite-dimensional subspace $E$ of $X$ and every
$\eps>0$, there is another finite-dimensional subspace $F$ of
$X$ containing $E$ such that the Banach-Mazur distance between
$F$ and $\ell_\infty^{(\dim(F))}$ is less than $1+\eps$. It is
well known \cite{Lazar-Lindenstrauss} that this property is
equivalent to the fact that $X^*$ is isometrically isomorphic
to an $L_1(\mu)$ space. The completions of the spaces
constructed in Theorem~\ref{str-conv-thm1} are ${\mathcal
L}^\infty_{1+}$-spaces, so they are preduals of $L_1(\mu)$
spaces. In the real case, it is known that preduals of
$L_1(\mu)$ spaces are almost-CL-spaces (see
\cite[Theorem~4.8]{Lin-book} or \cite[Corollary~3.6]{Lima2}).
Actually, the same is true for the complex case. We include
here a proof of this fact since we have been unable to find it
in the literature.

\begin{prop}
Let $X$ be a (real or complex) Banach space such that $X^*$ is
isometrically isomorphic to an $L_1(\mu)$ space. Then, $X$ is
an almost-CL-space.
\end{prop}

\begin{proof}
If we consider a maximal convex subset $F$ of $B_X$, the
Hahn-Banach and Krein-Milman theorems ensure that there is an
extreme point $f$ of the unit ball of $X^*=L_1(\mu)$ such that
$$
F=F(f):=\{x\in B_X\ : \ f(x)=1\}.
$$
We observe that the linear span of an extreme point $f$ in the
unit ball of an $L_1(\mu)$ space is an $L$-summand (i.e.\
$L_1(\mu)=\Lin(f)\oplus_1 Z$ for some closed subspace $Z$). So,
a result by \AA.~Lima \cite[Theorem~5.3]{Lima2} says that
$F(f)$ is not empty (we already knew it) and that $B_X$ is the
closure of $\aconv(F(f))$. This shows that $X$ is an
almost-CL-space.
\end{proof}

As a immediate consequence of this result we get the following.

\begin{corollary}
The completions of the non-complete lush strictly convex normed
spaces constructed in Theorem~\ref{str-conv-thm1} are
almost-CL-spaces and, therefore, they are not strictly convex.
\end{corollary}

We finish the section by remarking that the arguments of the
construction given in Theorem~\ref{str-conv-thm1} cannot be
adapted for smoothness, since a smooth norm on a subspace
cannot always be extended to the whole space (see
\cite[Theorem~8.3]{DGZ}).

\section{Separable lush spaces}\label{sec:lush=>stronglylush}
We have seen in the previous section that the completions of
the normed strictly convex lush spaces constructed are not
strictly convex by showing that they are almost-CL-spaces. We
cannot expect that every Banach space with numerical index $1$
is an almost-CL-space since there are lush spaces which do not
fulfil this property \cite[Example~3.4]{BKMW}. Nevertheless,
our aim here is to show that, in the separable case, lush
spaces actually have a much stronger property which in the real
case is very close to being an almost-CL-space and which will
allow us to show that a real lush space cannot be strictly
convex, unless it is one-dimensional.

We need a characterization of lushness given in
\cite{BKMM-lush} in terms of a norming subset of $S_{X^*}$.
Also, to carry some consequences to the non-separable case, we
need a result of the same paper saying that lushness is a
separably determined property. We state both results here for
easier reference.

\begin{prop}[\textrm{\cite[Theorems 4.1 and 4.2]{BKMM-lush}}] \label{complush}
Let $X$ be a Banach space.
\begin{enumerate}
\item[(a)] The following assertions are equivalent:
\begin{equi}
\item $X$ is lush\,.
\item  For every $x,y \in S_X$ and for every $\eps > 0$
    there is a slice $S=S(x^*,\eps) \subset B_X$, $x^* \in
    \ext[X]$,  such that
    $$x \in S \qquad \text{and} \qquad \dist
    \left(y,\aconv (S)\right) < \eps
    $$
    (i.e.\ $x^*$ in the definition of lushness can be
    chosen from $\ext[X]$).
\end{equi}
\item[(b)] The following two conditions are equivalent:
\begin{equi}
\item $X$ is lush,
\item Every separable subspace $E \subset X$ is contained in a
    separable lush subspace $Y$, $E \subset Y \subset X$.
\end{equi}
\end{enumerate}
\end{prop}

The following lemma, which will be the key to prove the main
result of the section, will be also useful in
section~\ref{sec:localizations} and does not depend upon the
separability of the space.

\begin{lemma}\label{lemma:lush-point-stronglylush}
Let $X$ be a lush space and let $K \subset B_{X^*}$ be the
weak* closure of $\ext[X^*]$ endowed with the weak* topology.
Then, for every $y\in S_X$, there is a $G_\delta$-dense subset
$K_y$ of $K$ such that $y\in \overline{\aconv (S(y^*,\eps))}$
for every $\eps > 0$ and every $y^* \in K_y$.
\end{lemma}

\begin{proof}
Fix $y \in S_X$. For every $n,m\in \N$, we consider
$$
K_{y,n,m} := \{x^* \in K\, :\, \dist \left(y,\aconv (S(x^*,1/n))\right) <
1/m\}.
$$
\indent \emph{Claim}.\ $K_{y,n,m}$ is weak*-open and dense in
$K$.

\noindent In fact, openness is almost evident: if $x^* \in
K_{y,n,m}$, then there is a finite set $A= \{a_1, \ldots a_k\}$
of elements of $S(x^*,1/n)$ such that $\dist \left(y,\aconv
(A)\right) < 1/m$. Denote
 $$
 U:= \{y^* \in K : \re y^*(a_i) > 1 -1/n \, \,{\textrm{ for  all }}\,\,
 i=1, \ldots, k
 \}.
 $$
$U$ is a weak*-neighborhood of $x^*$ in $K$, and $A \subset
S(y^*,1/n)$ for every $y^* \in U$. This means that $\dist
\left(y,\aconv (S(y^*,1/n))\right) < 1/m$ for all $y^* \in U$,
i.e.\ $U \subset K_{y, n,m}$.

To show density of $K_{y,n,m}$ in $K$, it is sufficient to
demonstrate that the weak* closure of $K_{y,n,m}$ contains
every extreme point $x^*$ of $S_{X^*}$. Since weak*-slices form
a base of neighborhoods of $x^*$ in $B_{X^*}$ (see
\cite[Lemma~3.40]{FHHMPZ}, for instance), it is sufficient to
prove that every weak*-slice $S(x,\delta)$, $\delta \in (0,
\min\{1/n,1/m\})$, intersects $K_{y,n,m}$, i.e.\ that there is
a point $y^* \in S(x,\delta) \cap K_{y,n,m}$. Which property of
$y^*$ do we need to make this true? We need that $y^*(x) > 1 -
\delta$, $y^* \in K$, and that $\dist \left(y,\aconv
(S(y^*,1/n))\right) < 1/m$. But the existence of such a $y^*$
is a simple application of item (ii) from
Proposition~\ref{complush}.(a). The claim is proved.

Now, we consider $K_y=\bigcap_{n,m\in\N} K_{y,n,m}$, which is a
weak*-dense $G_\delta$ subset of $K$ due to the Baire theorem.
\end{proof}

We are now ready to state and prove the main result of the
section.

\begin{theorem}\label{th:lush-stronglylush}
Let $X$ be a separable lush space. Then, there is a norming
subset $\tilde K$ of $S_{X^*}$ such that $B_{X} =
\overline{\aconv (S(x^*,\eps))}$ for every $\eps
> 0$ and for every $x^* \in \tilde K$. The last condition
implies that
$$
|x^{**}(x^*)|=1 \qquad \bigl(x^{**}\in\ext[X^{**}],\ x^*\in \tilde
K\bigr),
$$
and that in fact $\tilde K \subset \ext$.
\end{theorem}

\begin{proof}
We select a sequence $(y_n)$ dense in $S_X$ in such a way that
every element of the sequence is repeated infinitely many
times, and consider $\tilde K = \bigcap_{n \in \N} K_{y_n}$.
Due to the Baire theorem, $\tilde K$ is a weak*-dense
$G_\delta$ subset of $K$. This implies that for every $x \in
S_X$ and for every $\eps > 0$ there is an $x^* \in \tilde K$,
such that $x \in S(x^*,\eps)$ (i.e.\ $\tilde K$ is
$1$-norming). For $x_0^*\in \tilde K$ and $\eps>0$ fixed, the
inequality $\dist \left(y_n,\aconv (S(x_0^*,1/n))\right) < 1/n$
holds true for all $n \in \N$. Select an $N > 1/\eps$. Then,
for every $n > N$ we have $\dist \left(y_n,\aconv
(S(x_0^*,\eps))\right) < 1/n$. Since every element of the
sequence $(y_n)$ is repeated infinitely many times, this means
that $\dist \left(y_n,\aconv (S(x_0^*,\eps))\right) = 0$. So
the closure of $\aconv (S(x_0^*,\eps))$ contains the whole ball
$B_X$. Then,
$$
B_{X^{**}}=\overline{B_X}^{w^*} \subseteq
\overline{\aconv (S(x_0^*,\eps))}^{w^*}.
$$
Finally, the reversed Krein-Milman theorem gives us that
\begin{equation*}
\ext[X^{**}] \subset \overline{\T S(x_0^*,\eps)}^{w^*},
\end{equation*}
and the arbitrariness of $\eps>0$ gives us
\begin{equation*}
|x^{**}(x_0^*)|=1 \qquad \bigl(x^{**}\in\ext[X^{**}]\bigr).\qedhere
\end{equation*}
\end{proof}

We do not know whether the statement of the theorem is true in
the non-separable case.

\begin{remark}
From the proof of the above theorem it follows that the set
$\tilde K$ is actually a $G_\delta$-dense subset of the weak*
closure of $\ext[X^*]$ endowed with the weak* topology.
\end{remark}

As a consequence of the above theorem and results of \AA.~Lima
\cite{Lima2}, we get the following interesting version valid in
the real case.

\begin{corollary}\label{separable-lush-quasi-CL}
Let $X$ be a lush \textbf{real} separable space. Then, there is
a subset $A$ of $S_{X^*}$ norming for $X$ such that for every
$a^*\in A$ one has
$$
B_X=\overline{\aconv\bigl(\{x\in S_X\,:\, a^*(x)=1\}\bigr)}\,.
$$
\end{corollary}

\begin{proof}
By the above theorem, there is a subset $A$ of $S_{X^*}$
norming for $X$ such that
$$
|x^{**}(a^*)|=1 \qquad \bigl(x^{**}\in \ext[X^{**}],\ a^*\in A\bigr).
$$
Now, Theorems 3.1 and 3.5 of \cite{Lima2} give us that each
$a^*\in A$ attains its norm on $X$ and, moreover, that the
closed absolutely convex hull of the points of $B_X$ where
$a^*$ attains its norm is the whole ball $B_X$, as claimed.
\end{proof}

\begin{corollary}\label{cor-strictlyconvex-real-lush-onedim}
Let $X$ be a \textbf{real} Banach space which is lush. Then,
$X$ is neither strictly convex nor smooth, unless it is
one-dimensional.
\end{corollary}

\begin{proof}
If $X$ is a lush space, then every separable closed subspace
$Z$ of $X$ is contained in a separable lush subspace $Y$ by
Proposition~\ref{complush}.(b), and
Corollary~\ref{separable-lush-quasi-CL} provides us with a face
$F$ of $B_Y$ such that $B_Y=\overline{\conv(F\cup -F)}$. If $Y$
is not one-dimensional, then $F$ contains at least two distinct
points $y_1, y_2$ and $\frac12(y_1+y_2)\in F\subset S_Y$ is not
extreme. On the other hand, if $Y$ is not one-dimensional,
following the proof of \cite[Theorem~3.1]{Cheng-Li}, we get
that the smooth points of $F$ are exactly the norm-one elements
of the cone generated by $F$ which are not support points of
the cone. But then, the Bishop-Phelps theorem provides us with
(norm-one) support points of such a cone (see
\cite[Theorem~3.18]{Ph} for instance). Then, $F$ and so $S_Y$
contains non-smooth points. Therefore, $Y$ is not smooth, all
the more $X$.
\end{proof}

We do not know whether the above two results are true in the
complex case. We do not know either whether there are real
strictly convex Banach spaces with numerical index~$1$ others
than $\R$.

As a consequence of the corollary above, we get a negative
answer to a problem by M.~Popov, which he posed to the first
author in 1996 while discussing the still open problem on the
existence of a strictly convex Banach space with the Daugavet
property.

\begin{corollary}
A C-rich closed subspace of the \textbf{real} space $C[0,1]$ is
neither strictly convex nor smooth.
\end{corollary}

It is known that a subspace $X$ of $C[0,1]$ is C-rich whenever
$C[0,1]/X$ does not contain a copy of $C[0,1]$ (see
\cite[Proposition~1.2 and Definition~2.1]{KadPop}). Therefore,
the following is a particular case of the above proposition.

\begin{corollary}
Let $X$ be a closed subspace of the \textbf{real} space
$C[0,1]$. If $X$ is smooth or strictly convex, then $C[0,1]/X$
contains an isomorphic copy of $C[0,1]$.
\end{corollary}

Finally, another interesting consequence of
Theorem~\ref{th:lush-stronglylush} is the following.

\begin{corollary}
Let $X$ be an infinite-dimensional \textbf{real} Banach space
which is lush. Then $X^*$ contains an isomorphic copy of
$\ell_1$.
\end{corollary}

\begin{proof}
If $X$ is lush, by Proposition~\ref{complush}.(b), there is an
infinite-dimensional separable closed subspace $Y$ of $X$ which
is lush. Then, by Theorem~\ref{th:lush-stronglylush}, there is
a norming subset $\tilde K$ of $S_{Y^*}$ (in particular,
$\tilde K$ is infinite) such that
$$
|y^{**}(y^*)|=1 \qquad \bigl(y^{**}\in\ext[Y^{**}],\ y^*\in \tilde K\bigr).
$$
Now, Proposition~2 of \cite{LMP1999} shows that $Y^*$ contains
either $c_0$ or $\ell_1$. But a dual space contains
$\ell_\infty$ (hence also $\ell_1$) as soon as it contains
$c_0$ (see \cite[Theorem~V.10]{Dies} or
\cite[Proposition~2.e.8]{L-T}, for instance). Finally, if $Y^*$
contains a copy of $\ell_1$, then so does $X^*$ (see
\cite[p.~11]{vanDulst}, for instance).
\end{proof}

The above corollary has already been known for real spaces with
numerical index~$1$ which are Asplund or have the Radon-Nikod\'{y}m
property \cite{LMP1999}, and for real almost-CL-spaces
\cite{MartPaya-CL}.

\section{Extreme points of the unit
ball}\label{sec:localizations} The fact that a
finite-dimensional strictly convex Banach space with numerical
index~$1$ has to be one-dimensional is a direct consequence of
an old result by C.~McGregor \cite{Mc}. Namely, if $X$ is a
finite-dimensional Banach space with $n(X)=1$, then
$$
|x^*(x)|=1 \qquad \bigl(x^*\in \ext[X^*],\ x\in\ext[X] \bigr).
$$
In \cite[Problem~13]{KaMaPa} it was asked whether the above
result is also true in the infinite dimensional case. There are
two goals in this section. On the one hand, we will show that
this is not the case. We present two examples of lush spaces
(actually, C-rich subspaces of $C(K)$) such that there are
$x_0^*\in \ext[X^*]$ and $x_0\in \ext[X]$ with
$|x_0^*(x_0)|=0$. On the other hand, we will show that such an
example is not possible when the point $x_0$ is actually
extreme in $B_{X^{**}}$. In particular, we obtain that a lush
space which is MLUR or WMLUR must be one-dimensional.

Let us start with the first two examples. We give two different
constructions, one for both the real and the complex case and
another one for the complex case only, showing that the answer
to the already mentioned Problem~13 of \cite{KaMaPa} is
negative.

\begin{example}
{\slshape There is a C-rich subspace $X$ of the space $C[0,1]$
(hence $X$ is lush) and there are $f_0 \in \ext[X]$ and $x_0^*
\in \ext[X^*]$ satisfying $x_0^*(f_0)=0$.}
 \end{example}

\begin{proof}
Let us fix a function $f_{0} \in C[0,1]$ such that
$\|f_{0}\|=1$, $f_{0}(t)=0$ for $t \in (0, 1/3)$, and
$f_{0}(t)=1$ for $t \in (2/3, 1)$. We select a sequence of
intervals $\Delta_n \subset [0,1]$, $|\Delta_n| < 1/3$, such
that for every $(a,b) \subset [0,1]$ there is a $\Delta_j
\subset (a,b)$. Also, fix a null sequence $(\eps_n)$, $\eps_n
> 0$. Now one can easily construct functions
$f_n \in C[0,1]$, $n\in \N$, and functionals $f_n^* \in
C[0,1]^*$, $n = 0,1,2 \ldots$, recursively with the following
properties:
\begin{enumerate}
\item[(i)] $\|f_n\|=1$, $\|{f_n}|_{\Delta_n}\|_\infty = 1$
    and $\|f_n|_{[0,1] \setminus \Delta_n}\|_\infty \le
    \eps_n$.
\item[(ii)] All the $f_n$ are linear splines and $f_n(t)=0$
    in all the non-smoothness points of $f_k, \, k < n$, as
    well as in the points $0, 1/3, 2/3$ and $1$.
\item[(iii)] Every $f_n|_{(2/3, 1)}$ is linearly
    independent of $\{f_k|_{(2/3, 1)}\}_{k=0}^{n-1}$.
\item[(iv)] The measure representing $f_n^*$ is supported
    in $(2/3, 1)$.
\item[(v)] $f_n^*(f_m) = 0$ when $n \neq m$,  and
    $f_n^*(f_n) = 1$.
\end{enumerate}
Let us explain the construction. Since $f_{0}|_{(2/3, 1)} \neq
0$, we can select $f_0^*$ supported in $(2/3, 1)$ with
$f_0^*(f_0) =1$. Now we are going to select $f_1$. The
conditions (ii) and (v) on $f_1$ mean that the linear spline
$f_1$ must satisfy a finite number of linear equations:
$$
f_1(0) = f_1(1/3) = f_1(2/3) = 0, \, f_0^*(f_1) =0,
$$
so the set of splines supported on a fixed non-void interval
satisfying these conditions is a finite-codimensional subspace.
Select a norm-one spline $g_1\in C[0,1]$ supported on
$\Delta_1$ satisfying the equations above. Find a spline
$h_1\in C[0,1]$ of norm less than $\eps_1$, supported on $(2/3,
1) \setminus \Delta_1$, linearly independent of $f_0|_{(2/3, 1)
\setminus \Delta_1}$ and also satisfying the linear equations
for $f_1$. Then $f_1:= g_1 + h_1$ will serve its purpose. By
linear independence of $f_0|_{(2/3, 1) \setminus \Delta_1}$ and
$f_1|_{(2/3, 1) \setminus \Delta_1}= h_1$, we may find a
measure supported on $(2/3, 1) $ (and even more: on $(2/3, 1)
\setminus \Delta_1$) which annihilates $f_0$ and takes the
value $1$ on $f_1$. Take this measure as $f_1^*$. Then, in the
same way we construct $f_2$, then $f_2^*$, etc.

Now, we take $X:=\overline{\Lin}\{f_n\}_{n \in \N \cup \{0\}}$.
The property (i) ensures that $X$ is C-rich in $C[0,1]$. Thanks
to the property (ii), $\{f_n\}_{n \in \N \cup \{0\}}$ forms a
monotone basic sequence (i.e.\ a basis of $X$), and property
(v) gives us that the coordinate functionals can be written as
restrictions of $f_n^*$ to $X$.

Since $X$ is C-rich, there is a function $g \in S_X$ which
attains its norm only on $(0, 1/3)$. Fix $x_0^* \in \ext[X^*]$
with $x_0^*(g)=1$. Let $\mu \in S_{C[0,1]^*}$ be a measure
representing $x_0^*$. Then, $\mu$ is automatically supported on
$(0, 1/3)$, so $x_0^*(f_0)=0$. What remains to prove is that
$f_0 \in \ext[X]$. To this end, we consider $h \in X$ such that
$\|f_0 \pm h\| = 1.$ Since $f_{0}(t)=1$ for $t \in (2/3, 1)$,
we have that $h = 0$ on $(2/3, 1)$. But then, $h =
\sum_{n=0}^\infty f_n^*(h) f_n \equiv 0.$
\end{proof}

In the complex case, an easier example can be constructed.

\begin{example}
{\slshape There is a C-rich subspace $X$ of the complex space
$C(\T)$ (in particular, $X$ is lush and so $n(X)=1$), and
extreme points $x_0^*\in B_{X^*}$ and $\phi_0\in B_X$ such that
$x_0^*(\phi_0)=0$. }
\end{example}

\begin{proof}
Let $A(\D)$ be the disk algebra, considered as a closed
subspace of $C(\T)$. Then, $A(\D)$ is C-rich in $C(\T)$.
Indeed, let $\varphi(z)=\exp(z)/\e$ for every $z\in \D$. Then,
$\|\varphi\|=|\varphi(1)|=1$ and there is no other point on
$\T$ but $z=1$ where $\varphi$ attains its norm. Then, the
family
$$
\mathcal{A}=\{\varphi^n(z_0\,\cdot)\ : \ n\in\N,\ z_0\in \T\}
$$
belongs to $A(\D)$, and for every $\eps>0$ and every open
subset $U$ of $\T$, we may find an element of $\mathcal{A}$
which is at $\eps$-distance from a function whose support is
inside $U$.

Let $X=\Lin \bigl\{A(\D),\ \phi_0 \bigr\}$, where $\phi_0$ is
any function in $C(\T)$ for which there are non-empty open sets
$U_1$, $U_2$ and $U_3$ of $\T$ such that $\phi_0\equiv 1$ on
$U_1$, $\phi_0\equiv -1$ on $U_2$ and $\phi_0\equiv 0$ on
$U_3$. Then, $X$ is C-rich because it contains $A(\D)$. Next,
$\phi_0$ is extreme on $B_X$. Indeed, if $g=\alpha \phi_0 +
f\in X$ is such that $\|\phi_0 \pm g\|\leq 1$, then $g\equiv 0$
on $U_1\cup U_2$, so $f\equiv \alpha$ on $U_1$ and $f\equiv
-\alpha$ on $U_2$. It follows that $\alpha=0$ and so $g=0$.
Also, for every $z\in \T$, the functional $\delta_z$ is extreme
in $B_{X^*}$. Namely, for every $z\in \T$ there is a function
$\varphi \in A(\D)$ which attains its norm only at the point
$z$. Then, there is an extreme point $x^*$ of $B_{X^*}$ such
that $|x^*(\varphi)| = 1$. Then, the norm-one measure which
represents $x^*$ must be supported on $z$ (otherwise the
integral would be strictly smaller than $1$), so $x^*$ is of
the form $\theta \delta_z$. Finally, taking $z\in U_2$ and
calling $x_0^*=\delta_z$, we have that $x_0^*\in \ext[X^*]$ and
$x_0^*(\phi_0)=0$.
\end{proof}

Let us comment that the extreme points $f_0$ and $\phi_0$ of
the examples above are not rotund. In fact, we do not know if a
rotund point may exist in a lush space with dimension greater
than one. We recall that a point $x$ in the unit sphere of a
Banach space $X$ is said to be \emph{rotund} if it is not an
element of any nontrivial closed segment in the unit sphere or,
equivalently, if $\|x+y\|=2$ for some $y\in B_X$ implies $y=x$.

To finish the section, we show that in the previous examples
the extreme points of the unit ball cannot be $w^*$-extreme. We
will use this to show that there are no lush spaces which are
WMLUR.

\begin{prop}\label{prop-w*extreme}
Let $X$ be a lush space. Then, for every $w^*$-extreme point
$x_0$ of $B_X$ and every $x^*\in \ext[X^*]$, one has
$|x^*(x_0)|=1$. In particular, this happens for WMLUR points of
$B_X$.
\end{prop}

\begin{proof}
We fix a $w^*$-extreme point $x_0$ of $B_X$. By
Lemma~\ref{lemma:lush-point-stronglylush}, there is a subset
$K_{x_0}$ of $\ext[X^*]$ norming for $X$ such that $x_0 \in
\overline{\aconv (S(k^*,\eps))}$ for every $\eps > 0$ and every
$k^* \in K_{x_0}$. Then, since $x_0$ is an extreme point of the
bidual ball, the same argument as at the end of the proof of
Theorem~\ref{th:lush-stronglylush} shows that for all $k^*\in
K_{x_0}$
$$
|x_0(k^*)|=1.
$$
Since $K_{x_0}$ is norming for $X$, we have that $B_{X^*}$ is
the weak$^*$-closure of $\aconv(\tilde{K})$, and the reversed
Krein-Milman theorem gives us that the set $\ext[X^*]$ is
contained in the $w^*$-closure of $\T K_{x_0}$. The result
follows since $x_0\in X$.
\end{proof}

As a consequence, we have the following prohibitive result. In
the real case, it is a particular case of
Corollary~\ref{cor-strictlyconvex-real-lush-onedim}, since
WMLUR spaces are strictly convex.

\begin{corollary}
Let $X$ be a lush space. Then, $X$ is not WMLUR (in particular,
it is not MLUR), unless it is one-dimensional.
\end{corollary}

\begin{proof}
Since $X$ is a WMLUR, Proposition~\ref{prop-w*extreme} gives us
that $|x^*(x)|=1$ for every $x^*\in \ext[X^*]$ and every $x\in
S_X$. But this clearly implies that $X$ is one-dimensional.
\end{proof}

Let us mention that another consequence of
Proposition~\ref{prop-w*extreme} is that every $w^*$-extreme
point of a lush space is actually MLUR, as the following remark
shows.

\begin{remark}
{\slshape Let $X$ be a Banach space and let $x$ be a point in
$B_X$ so that $|x^*(x)|=1$ for every $x^*\in \ext[X^*]$. Then,
$x$ is an MLUR point of $B_X$. Indeed, fixed $y\in X$, we take
$x^*\in \ext[X^*]$ so that $x^*(y)=\|y\|$ and we estimate as
follows
\begin{equation*}
\underset{\pm}{\max}\|x\pm y\|\geq\underset{\pm}{\max}\big|x^*(x)\pm
\|y\|\big|\geq \big(|x^*(x)|^2+\|y\|^2\big)^{1/2}=(1+\|y\|^2)^{1/2}\,.
\end{equation*} }
\end{remark}

Finally, we present an example showing that the results above
are not valid for Banach spaces having the alternative Daugavet
property. We do not know whether they are true for spaces with
numerical index $1$.

\begin{example}
{\slshape The real or complex space $X=C\big([0,1],
\ell_2^2\big)$ has the alternative Daugavet property (and even
the Daugavet property). However, there exist $x_0^*\in
\ext[X^*]$ and a MLUR point $f_0$ of $B_X$ such that
$|x_0^*(f_0)|<1$.}
\end{example}

\begin{proof}
First, $C([0,1],\ell_2^2)$ has the alternative Daugavet
property by \cite[Theorem~3.4]{MaOi}, for instance. Now, we fix
any $x_0\in S_{\ell_2^2}$ and consider $f_0\in S_X$ given by
$f_0(t)=x_0$ for every $t\in [0,1]$. To prove that $f_0$ is an
MLUR point in $B_X$, we take $g\in X$ and we observe that
\begin{align*}
\underset{\pm}{\max}\|f_0\pm
g\|=\underset{t\in[0,1]}{\sup}\underset{\pm}{\max}\|x_0\pm g(t)\|\geq
\underset{t\in[0,1]}{\sup} \bigl(1+\|g(t)\|^2\bigr)^{1/2}=
(1+\|g\|^2)^{1/2}.
\end{align*}
We notice that the above inequality becomes an equality when
one considers $g\in X$ given by $g(t)=x_0^\perp$ for every
$t\in[0,1]$, where $x_0^\perp\in S_{\ell_2^2}$ is orthogonal to
$x_0$. Finally, it suffices to take $x_0^*\in \ext[X^*]$ so
that $x_0^*(g)=1$ to get the desired condition. Indeed,
\begin{equation*}
\sqrt{2}=\underset{\omega\in\T}{\max}\|f_0+ \omega\,g\|\geq
\underset{\omega\in\T}{\max}|x_0^*(f_0)+\omega| =1+|x_0^*(f_0)|,
\end{equation*}
so $|x_0^*(f_0)|\leq \sqrt{2}-1<1$.
\end{proof}

\emph{Acknowledgement.} The authors thank the anonymous referee
for multiple stylistical improvements.

\end{document}